\documentclass[12pt]{amsart}
\usepackage{amssymb}
\usepackage[matrix,arrow,curve]{xy}
\newtheorem{Theorem}{Theorem}[section]
\newtheorem{Proposition}[Theorem]{Proposition}
\newtheorem{Lemma}[Theorem]{Lemma}
\newtheorem{Definition}[Theorem]{Definition}
\newtheorem{Remark}[Theorem]{Remark}
\newtheorem{Claim}{Claim}
\title[An isoperimetric inequality on the projective plane]{An isoperimetric 
inequality for 
the second non-zero eigenvalue of the Laplacian
on the projective plane}
\author{Nikolai S. Nadirashvili}
\address[Nikolai S. Nadirashvili]{CNRS, I2M UMR 7353 --- 
Centre de Math\'e\-ma\-ti\-ques et Informatique,
Marseille, France}
\email{nikolay.nadirashvili@univ-amu.fr}
\author{Alexei V. Penskoi}\thanks{The work of the second author
was partially supported by the Simons Foundation and 
by the Young Russian Mathematics award.}
\address[Alexei V. Penskoi]{Department of Higher Geometry and Topology, 
Faculty of Mathematics and Mechanics, Moscow State University,
Leninskie Gory, GSP-1, 119991, Moscow, Russia \newline {\em and}
\newline Faculty of Mathematics,
National Research University Higher School of Economics,
6 Usacheva Str., 119048, Moscow, Russia \newline {\em and}
\newline Interdisciplinary Scientific Center
J.-V. Poncelet (ISCP, UMI 2615), Bolshoy Vlasyevskiy 
Pereulok 11, 119002, Moscow, Russia}
\email[Corresponding author]{penskoi@mccme.ru}
\date{}
\DeclareMathOperator{\ord}{ord}
\DeclareMathOperator{\sgn}{sgn}
\DeclareMathOperator{\Area}{Area}
\DeclareMathOperator{\tr}{tr}
\begin{document}
\begin{abstract}
We prove an isoperimetric inequality for the
second non-zero eigenvalue of the Laplace-Beltrami 
operator on the real projective plane. For a
metric of unit area this eigenvalue is not greater than~$20\pi.$
This value is attained in the limit by a sequence of metrics of area one 
on the projective plane. The limiting metric is singular and could be 
realized as a union of the projective plane and the sphere touching at a 
point, with standard metrics and the ratio of the areas~$3:2.$
It is also proven that the multiplicity of the second non-zero eigenvalue 
on the projective plane is at most~$6.$
\end{abstract}
\maketitle
\section{Introduction}
Let $M$ be a closed surface and $g$ be
a Riemannian metric on $M.$
Let us consider the Laplace-Beltrami
operator $\Delta:C^\infty(M)\longrightarrow C^\infty(M)$
associated with the metric $g,$
$$
\Delta f=-\frac{1}{\sqrt{|g|}}\frac{\partial}{\partial x^i}%
\left(\sqrt{|g|}g^{ij}\frac{\partial f}{\partial x^j}\right),
$$
and its eigenvalues
\begin{equation}\label{eigenvalues}
0=\lambda_0(M,g)<\lambda_1(M,g)\leqslant%
\lambda_2(M,g)\leqslant\lambda_3(M,g)\leqslant\dots
\end{equation}
Let us denote by $m(M,g,\lambda_i)$
the multiplicity of the eigenvalue $\lambda_i(M,g),$
i.e. how many times the value of $\lambda_i(M,g)$
appears in the sequence~\eqref{eigenvalues}.

Let us consider a functional
$$
\bar{\lambda}_i(M,g)=\lambda_i(M,g)\Area(M,g),
$$
where $\Area(M,g)$ is the area of $M$ with respect to the
Riemannian metric $g.$ This functional is sometimes called
an eigenvalue normalized by the area or simply a normalized eigenvalue.

Yang and Yau proved in the paper~\cite{Yang-Yau1980} that
if $M$ is an orientable surface of genus $\gamma$ then
$$
\bar{\lambda}_1(M,g)\leqslant 8\pi(\gamma+1).
$$
Actually, the arguments of Yang and Yau imply a stronger estimate, 
$$
\bar{\lambda}_1(M,g)\leqslant 8\pi\left[\frac{\gamma+3}{2}\right],
$$
see the paper~\cite{ElSoufi-Ilias1984}
and also~\cite{Nadirashvili1996}.
Here $[\cdot]$ denotes the integer part of a number.

Later Korevaar proved in the paper~\cite{Korevaar1993}
that there exists a constant $C,$ such that for any 
$i>0$ and any compact surface $M$ of genus $\gamma$
the following upper bound holds:
$$
\bar{\lambda}_i(M,g)\leqslant C(\gamma+1)i.
$$
Recently this upper bound was improved by Hassannezhad~\cite{Hassannezhad2011}. 
She proved that there exists a constant $C,$ such that for any 
$i>0$ and any compact surface $M$ of genus $\gamma$,
the following upper bound holds:
$$
\bar{\lambda}_i(M,g)\leqslant C(\gamma+i).
$$
It follows that the functionals $\bar{\lambda}_i(M,g)$
are bounded from above and it is a natural question to find
for a given compact surface $M$ and number $i\in\mathbb{N}$
the quantity
$$
\Lambda_i(M)=\sup_g\bar\lambda_i(M,g),
$$
where the supremum is taken over the space of all
Riemannian metrics $g$ on $M.$

Let us remark that the functional $\bar{\lambda}_i(M,g)$
is invariant under rescaling of the metric 
$g\mapsto tg,$ where $t\in\mathbb{R}_+.$ It follows that it is
equivalent to the problem of finding $\sup\lambda_i(M,g),$
where the supremum is taken over the space of all
Riemannian metrics $g$ of area $1$ on $M.$
That's why this problem is sometimes called the {\em isoperimetric
problem for eigenvalues of the Laplace-Beltrami operator.}

\begin{Definition}\label{maximal-closed}
Let $M$ be a closed surface.
A metric $g_0$ on $M$ is called maximal for the functional
$\bar{\lambda}_i(M,g)$ if
$$
\Lambda_i(M)=\bar{\lambda}_i(M,g_0)
$$
\end{Definition}

If a maximal metric exists, it is defined up to
multiplication by a positive constant due to the rescaling
invariance of the functional.

Surprisingly, the list of known results is quite short.

Hersch proved in 1970 in the paper~\cite{Hersch1970}
that the standard metric on the sphere
is the unique maximal metric for
$\bar{\lambda}_1(\mathbb{S}^2,g)$
and 
$$
\Lambda_1(\mathbb{S}^2)=8\pi.
$$

Li and Yau proved in 1982 in the paper~\cite{Li-Yau1982}
that the standard metric on the projective plane
is the unique maximal metric for
$\bar{\lambda}_1(\mathbb{RP}^2,g)$
and 
$$
\Lambda_1(\mathbb{RP}^2)=12\pi.
$$

The first author proved in 1996 in the paper~\cite{Nadirashvili1996}
that the standard metric on the equilateral torus
is the unique maximal metric for
$\bar{\lambda}_1(\mathbb{T}^2,g)$
and 
$$
\Lambda_1(\mathbb{T}^2)=\frac{8\pi^2}{\sqrt{3}}.
$$

It is not always that a maximal metric exists. As it was proved by the first
author in 2002 in the paper~\cite{Nadirashvili2002}
and later with a different argument by Petrides~\cite{Petrides2014},
$$
\Lambda_2(\mathbb{S}^2)=16\pi.
$$
However, there is no maximal metric. The supremum is attained
as a limit on a sequence of smooth metrics on the sphere converging
to a singular metric on two spheres of the same radius touching
in a point.

The functional $\bar{\lambda}_i(M,g)$ depends continuously
on the metric $g.$ 
However, when $\bar{\lambda}_i(M,g)$ is a multiple eigenvalue 
this functional is not in general differentiable. If we consider an
analytic variation $g_t$ of the metric $g=g_0,$
then it was proved by Berger~\cite{Berger1973},
Bando and Urakawa~\cite{Bando-Urakawa1983},
El Soufi and Ilias~\cite{ElSoufi-Ilias2008}
that the left and right derivatives
of the functional $\bar{\lambda}_i(M,g_t)$ with respect to $t$ exist.
This leads us to the  following definition
given by the first author in the paper~\cite{Nadirashvili1996}
and by El Soufi and Ilias in the 
papers~\cite{ElSoufi-Ilias2000,ElSoufi-Ilias2008}.

\begin{Definition}
A Riemannian metric $g$ on  a closed surface
$M$ is called extremal metric for the
functional $\bar\lambda_i(M,g)$ if for any analytic deformation
$g_t$ such that $g_0=g$ one has
$$
\frac{d}{dt}\bar{\lambda}_i(M,g_t)%
\left.\vphantom{\raisebox{-0.5em}{.}}\right|_{t=0+}\leqslant0%
\leqslant\frac{d}{dt}\bar{\lambda}_i(M,g_t)%
\left.\vphantom{\raisebox{-0.5em}{.}}\right|_{t=0-}.
$$
\end{Definition}

Jakobson, the first author and Polterovich
proved in 2006 in the paper~\cite{Jakobson-Nadirashvili-Polterovich2006}
that the metric on the Klein bottle realized as so called
bipolar Lawson surface $\tilde{\tau}_{3,1},$
is extremal for $\bar\lambda_1(\mathbb{KL},g).$
It was conjectured in this paper that this metric is 
the maximal one. 
El Soufi, Giacomini and Jazar proved in the same year
in the paper~\cite{ElSoufi-Giacomini-Jazar2006}
that this metric on $\tilde{\tau}_{3,1}$
is the unique extremal metric for the
$\bar\lambda_1(\mathbb{KL},g).$ 
It follows from the 
results of~\cite{Matthiesen-Siffert2017q} 
that there exists a  smooth (up to at most a finite 
number of conical points)  metric $g_K$ on the Klein bottle
such that $\sup\bar{\lambda}_1(\mathbb{KL},g)$ is attained on $g_K.$
It could be then shown (a detailed exposition of this argument 
could be found in~\cite{Cianci-Medvedev2017q})
that the metric on $\tilde{\tau}_{3,1}$ is the maximal one
and, hence,
$$
\Lambda_1(\mathbb{KL})=\bar{\lambda}_1(\mathbb{KL},g_{\tilde{\tau}_{3,1}})=%
12\pi E\left(\frac{2\sqrt{2}}{3}\right),
$$
where $E$ is the complete elliptic integral of the second kind
and $g_{\tilde{\tau}_{3,1}}$ is the metric on $\tilde{\tau}_{3,1}.$

More results on extremal metrics on tori and Klein bottles
could be found in the 
papers~\cite{ElSoufi-Ilias2000,Karpukhin2013,Karpukhin2014,
Karpukhin2015,Lapointe2008,Penskoi2012,Penskoi2013,Penskoi2015}. 
A review of these
results is given by the second author in the paper~\cite{Penskoi2013a}.

It was shown in \cite{JLNNP} using a combination of analytic and numerical tools
that the maximal metric for the first eigenvalue on 
the surface of genus two is the metric on the Bolza surface $\mathcal P$
induced from the canonical metric on the sphere
using the standard covering ${\mathcal P}\longrightarrow\mathbb{S}^2.$ 
The authors stated this result as a conjecture, because
the argument is partly  based on a numerical calculation. The proof of
this conjecture was given in a recent preprint~\cite{Nayatani-Shoda2017q}.

The first author and Sire proved in 2015 in the 
paper~\cite{Nadirashvili-Sire2015q2} the equality
$$
\Lambda_3(\mathbb{S}^2)=24\pi.
$$
It turns out that there is no maximal metric
but the supremum could be
obtained as a limit on a sequence of metrics on the sphere
converging
to a singular metric on three touching spheres of the same radius.
It was conjectured in the paper~\cite{Nadirashvili2002,Nadirashvili-Sire2015q2}
that
$$
\Lambda_k(\mathbb{S}^2)=8\pi k.
$$
This conjecture was proven in the recent 
paper~\cite{Karpukhin-Nadirashvili-Penskoi-Polterovich2017q}
by Karpukhin, Polterovich and the authors.

The main goal of the present paper is to prove the
following result.

\begin{Theorem}\label{maintheorem}
The supremum of the normalized second nonzero eigenvalue on the projective 
plane over the space of all Riemannian metrics on $\mathbb{RP}^2$ is given 
by 
\begin{equation}\label{mainequation}
\Lambda_2(\mathbb{RP}^2)=20\pi.
\end{equation}
There is no maximal metric, even among metrics with a finite number of 
conical singularities. The supremum is attained 
in the limit by a sequence of metrics of area one 
on the projective plane. The limiting metric is singular and could be 
realized as a union of the projective plane and the sphere touching at a 
point, with standard metrics and the ratio of the areas~$3:2.$
\end{Theorem}

We postpone the definition of metrics with conical singularities
till Section~\ref{harmonic-conical}.

\begin{Remark}
This Theorem could be stated as an isoperimetric inequality
$$
\lambda_2(\mathbb{RP}^2,g)\leqslant20\pi
$$
for any metric $g$ of area $1.$
\end{Remark}

\begin{Remark}
It would be interesting to check whether the equality in~\eqref{mainequation} 
could be attained in the limit only by a sequence of metrics converging to  
a union of touching projective plane and sphere
with standard metrics and the ratio of the areas~$3:2$,
or there exist other maximizing sequences.
\end{Remark}

\begin{Remark}
The degenerating sequence of metrics in Theorem \ref{maintheorem} 
illustrates the "bubbling phenomenon'' arising in the maximization 
of higher eigenvalues, see \cite{Nadirashvili-Sire2015b} for details.
\end{Remark}

\begin{Remark}
It was conjectured in the 
paper~\cite{Karpukhin-Nadirashvili-Penskoi-Polterovich2017q},
written after the first version of the present paper,
that the equality
$$
\Lambda_k(\mathbb{RP}^2)=4\pi (2k+1)
$$
holds for any $k\geqslant 1$.
\end{Remark}

The paper is organized in the following way.
In Section~\ref{minimal_immersions} we recall the relation
between extremal metrics and minimal immersions into spheres
and explain the importance of upper bounds on the multiplicities
of eigenvalues.

In Section~\ref{nodal} we recall the basics of the theory of nodal
graphs and the Courant Nodal Domain Theorem. We use them 
in Section~\ref{multiplicity} in order to obtain  an upper bound
for the multiplicity $m(\mathbb{RP}^2,g,\lambda_2).$
Let us remark that bounds on multiplicity of eigenvalues
of the Laplace-Beltrami operator on surfaces
were subject of numerous papers, see 
e.g.~\cite{Besson1980,Cheng1976,
Hoffmann-Ostenhof-Hoffmann-Ostenhof-Nadirashvili1999,
Karpukhin-Kokarev-Polterovich2014,Nadirashvili1987}.

In Section~\ref{harmonic-conical} we pass from minimal immersions
to harmonic immersions, extend our considerations to harmonic
immersions with branch points
and metrics with conical singularities and explain why the results
from the previous sections also hold in this case.

In Section~\ref{Calabi-Barbosa-implications} we recall the Calabi-Barbosa
Theorem about harmonic immersions with branch points
$\mathbb{S}^2\longrightarrow\mathbb{S}^n$
and apply it to our situation.

Section~\ref{space_of_harmonic_maps} contains the description
of the space of harmonic immersions with branch points
$\mathbb{S}^2\longrightarrow\mathbb{S}^4$ due to Bryant
and results about singularities of these maps.

Section~\ref{existence} deals with the question
of existence of maximal metrics.

Finally, in Section~\ref{mainproof} we prove Theorem~\ref{maintheorem}.

\noindent {\bf Acknowledgements.} The authors are very indebted to Mikhail
Kar\-pu\-khin, Iosif Polterovich and the referee for useful remarks and
suggestions.
The second author is very grateful to
the Institut de Ma\-th\'e\-ma\-ti\-ques de Marseille (I2M, UMR~7373)
for the hospitality. The second author is very indebted to Pavel Winternitz
for fruitful discussions.

\section{Extremal metrics and minimal 
immersions into spheres}\label{minimal_immersions}

In this Section we recall the relation
between extremal metrics and minimal immersions into spheres
and explain the importance of upper bounds on the multiplicities
of eigenvalues.

Let us recall the definition of a minimal map, see 
e.g.~\cite{Eells-Lemaire1978,Eells-Sampson1964}.

Let $(M,g)$ be a Riemannian manifold of dimension $m.$
Let $\alpha$ be a symmetric bilinear $2$-form
on $TM.$ Let $\sigma_k$ be the $k$-th elementary
symmetric function.
Let $\sigma_k(\alpha)=\sigma_k(\lambda_1,\dots,\lambda_m),$
where $\lambda_i$ are eigenvalues of $\alpha$ related to $g,$
i.e. the roots of the polynomial $\det(\alpha_{ij}-\lambda g_{ij})=0.$

\begin{Definition}
Let $(M,g)$ and $(N,h)$ be Riemannian manifolds.
A smooth
map $f:M\longrightarrow N$ is called
minimal if $f$ is an extremal for the volume functional
$$
V[f]=\int_M\sqrt{|\sigma_m(f^*h)|}\,dVol_g,
$$
where $m=\dim M.$  
\end{Definition}

It is well-known that a surface 
$M\looparrowright\mathbb{R}^3$ is minimal if and only if
the coordinate functions $x^i$ are harmonic
with respect to the Laplace-Beltrami operator on $M.$
A similar result holds for a minimal submanifold in
$\mathbb{R}^n.$ Since harmonic functions are eigenfunctions
with eigenvalue~$0,$ it is natural to ask what is an analogue
of this statement for a non-zero eigenvalue. The answer was
given by Takahashi in 1966.

\begin{Theorem}[Takahashi~\cite{Takahashi1966}]\label{takahashi}
If an isometric immersion 
$$
f:M\looparrowright\mathbb{R}^{n+1},\quad f=(f^1,\dots,f^{n+1}),
$$
is defined by eigenfunctions
$f^i$ of the Laplace-Beltrami operator $\Delta$
with a common eigenvalue $\lambda,$
$$
\Delta f^i=\lambda f^i,
$$
then (i) the image $f(M)$ lies on the sphere $\mathbb{S}^n_R$
of radius $R$ with the center at the origin
such that 
\begin{equation}\label{lambda-R}
\lambda=\frac{\dim M}{R^2},
\end{equation}
(ii) the immersion $f:M\looparrowright\mathbb{S}^n_R$
is minimal.

If 
$$
f:M\looparrowright\mathbb{S}^{n}_R\subset\mathbb{R}^{n+1},
\quad f=(f^1,\dots,f^{n+1}),
$$ 
is a minimal isometric immersion of a manifold $M$
into the sphere $\mathbb{S}^{n}_R$
of radius $R,$ then $f^i$
are eigenfunctions of the Laplace-Beltrami
operator~$\Delta,$
$$
\Delta f^i=\lambda f^i,
$$
with the same eigenvalue $\lambda$
given by formula~\eqref{lambda-R}.
\end{Theorem}

We introduce the eigenvalue counting function
$$
N(\lambda)=\#\{\lambda_i|\lambda_i<\lambda\}.
$$
Takahashi's Theorem~\ref{takahashi} implies that if $M$
is isometrically minimally immersed in the sphere $\mathbb{S}^{n}_R,$
then among the eigenvalues of $M$ there are at least $n+1$
eigenvalues equal to $\frac{\dim M}{R^2}.$ It is easy
to see that $\lambda_{N\left(\frac{\dim M}{R^2}\right)}$ is the first
eigenvalue equal to $\frac{\dim M}{R^2}.$ This is important due to the
following theorem.

\begin{Theorem}[Nadirashvili~\cite{Nadirashvili1996},
El Soufi, Ilias~\cite{ElSoufi-Ilias2008}]\label{Nadirashvili-ElSoufi-Ilias}
Let $M\looparrowright\mathbb{S}^n_R$ be an immersed minimal
compact submanifold. Then the metric induced on $M$ by this immersion
is extremal for the functional
$\bar{\lambda}_{N\left(\frac{\dim M}{R^2}\right)}(M,g).$

If a metric on a compact manifold $M$ is extremal for
some eigenvalue then there exists
an isometric minimal immersion to the sphere
$M\looparrowright\mathbb{S}^n_R$
by eigenfunctions with eigenvalue
$\frac{\dim M}{R^2}$ of the Laplace-Beltrami operator corresponding
to this metric.
\end{Theorem}

If a metric is extremal for $\bar{\lambda}_i(M,g),$
then there exist a minimal immersion of $M$ by corresponding eigenfunctions
into $\mathbb{S}^n\subset\mathbb{R}^{n+1}.$ If the image is
not contained in some hyperplane then one should have at least
$n+1$ linearly indepenent eigenfunctions. This means that
$n+1\leqslant m(M,g,\lambda_i).$

If follows that if we have an upper bound on the multiplicity
of an eigenvalue then we have an upper bound on the dimension
of the sphere where $M$ is minimally immersed by the
corresponding eigenfunctions.

We later use Theorem~\ref{Nadirashvili-ElSoufi-Ilias} for 
$M=\mathbb{RP}^2.$ In this case $\dim M=2.$ Using rescaling one can consider only
the case of $R=1.$ Remark that Theorem~\ref{Nadirashvili-ElSoufi-Ilias}
holds also for a non-orientable $M.$

Since we are interested in the functional $\bar{\lambda}_2(\mathbb{RP}^2,g),$
we need an upper bound for $m(\mathbb{RP}^2,g,\lambda_2)$
in order to bound the dimension of the sphere which is sufficient
to consider.

\section{Nodal graphs and Courant Nodal Domain Theorem}\label{nodal}

In this Section we recall the basics of the theory of nodal
graphs and the Courant Nodal Domain Theorem that we need in order
to obtain in Section~\ref{multiplicity} an upper bound
$m(\mathbb{RP}^2,g,\lambda_2)\leqslant 6.$

Let us now recall the following theorem due to Bers.

\begin{Theorem}[L. Bers \cite{Bers1955}]\label{BersTheorem}
Let $(M,g)$ be a compact 2-dimensional closed Riemannian
manifold and $x_0$ is a point on $M.$ Then there
exist its neighbourhood chart $U$ 
with coordinates $x=(x^1,x^2)\in U\subset\mathbb{R}^2$
centered at $x_0$ such that for any 
eigenfunction $u$ of the Laplace-Beltrami operator on $M$
there exists an integer $n\geqslant0$ and a non-trivial homogeneous harmonic 
polynomial $P_n(x)$ of degree $n$ on the Euclidean plane $\mathbb{R}^2$ such 
that 
$$
u(x)=P_n(x)+O(|x|^{n+1}),
$$
where $x\in U$.
\end{Theorem}

The integer $n$ from Bers's Theorem~\ref{BersTheorem}
is called an order of zero of an eigenfunction $u$ at a point $x_0.$
Let us denote it by $\ord_{x_0}u.$

Consider the sets 
$$
N^l(u)=\{x\in M|\ord_x u\geqslant l\}.
$$

\begin{Definition}
The set $N^1(u)$ is called a nodal set of $u.$
Connected components of its complement $M\setminus N^1(u)$
are called the nodal domains of~$u.$
\end{Definition}

It is well-known that in the polar 
coordinates $r,\varphi$ in $\mathbb{R}^2$ any homogeneous
harmonic polynomial $P_n$ of degree $n$ has the form 
\begin{equation}\label{harmonic_pol}
P_n(r,\varphi)=r^n(A\cos n\varphi+B\sin n\varphi).
\end{equation}
The zeroes of such polynomials form $n$ 
straight lines intersecting at origin at equal angles.

It follows that the nodal set $N^1(u)$ is a graph such that the points of
$N^2(u)$ are its vertices and the
connected components of $N^1(u)\setminus N^2(u)$
are its edges. 
\begin{Definition}
This graph is called a nodal graph of an 
eigenfunction~$u.$
\end{Definition}

Let us remark that if $x_0$ is a vertex of the nodal graph then
it is a zero of $u$ and there is $2\ord_{x_0}u$ edges emanating from $x_0$
in a sufficiently small neighborhood of $x_0.$ Globally some of these edges could
form loops starting and ending at $x_0.$

Let us remark that {\em locally in a neighborhood of zero $x_0$ of order $n$
the nodal graph $N^1(u)$ looks like a star consisting of $2n$ rays
with equal angles between adjacent rays}. 
Let us give the following definition in order to be more precise.
\begin{Definition}
A star $S_{x_0}(N^1(u))$ at the vertex $x_0$ of the nodal graph
$N^1(u)$ of an eigenfunction $u$ consists of $2n$ unitary tangent vectors to 
edges emanating from $x_0,$ where $n$ is the order of zero of $u$
at $x_0.$
\end{Definition}

It follows from formula~\eqref{harmonic_pol} 
that in coordinates given by the Bers Theorem~\ref{BersTheorem}
the angles between adjacent vectors in $S_{x_0}(N^1(u))$ are equal.

If one has a triangulation of a surface $M$ with $V$ vertices,
$E$ edges and $F$ faces, then one has the well-known formula
for the Euler characteristic,
\begin{equation}\label{euler-equation}
\chi(M)=V-E+F.
\end{equation}

Let us consider an eigenfunction $u.$
If we consider the vertices of a nodal graph, 
the edges of a nodal graph and the nodal domains of a function $u,$
then the formula~\eqref{euler-equation} does not in general
hold since the nodal domains are not in general homeomorphic
to a disc. As a result, we obtain in this case
only the Euler inequality
\begin{equation}\label{euler-inequality}
\chi(M)\leqslant V-E+F
\end{equation}
that implies the following well-known Lemma.
\begin{Lemma}\label{nodal-inequality}
Let $u$ be an eigenfunction.
Let $x_j,$ $j=1,\dots,n,$
be zeroes of $u$ of order $m_j>1.$
Let $\Omega_j,$ $j=1,\dots,s,$
be nodal domains of the function $u.$
Then
$$
s\geqslant\chi(M)-n+\sum\limits_{j=1}^nm_j.
$$
\end{Lemma}

\noindent{\bf Proof.} One can immediately see 
that $V=n,$ $F=s.$
Since $2\ord_{x_j}u=2m_j$ edges emanate from $x_j$
and each edge connects two vertices,
one has $E=\sum\limits_{j=1}^{n}m_j.$
It is sufficient now to apply
inequality~\eqref{euler-inequality}.
$\Box$

Let us now recall the following theorem 
(remark that we count eigenvalues starting
from $\lambda_0$).

\begin{Theorem}[Courant Nodal Domain Theorem~\cite{Courant-Hilbert1931}]\label{CNDT}
An eigenfunction corresponding to the eigenvalue $\lambda_i$ has at most
$i+1$ nodal domains.
\end{Theorem}

Lemma~\ref{nodal-inequality} and Courant Nodal Domain Theorem~\ref{CNDT}
imply immediately the following Proposition.

\begin{Proposition}
Let $u$ be an eigenfunction corresponding 
to the eigenvalue $\lambda_i.$
Let $x_j,$ $j=1,\dots,n,$
be zeroes of $u$ of order $m_j>1.$
Then
\begin{equation}\label{nodal-inequality-1}
i+1\geqslant\chi(M)-n+\sum\limits_{j=1}^nm_j.
\end{equation}
\end{Proposition}

\section{Multiplicity of the second non-zero
eigenvalue of the Laplace-Beltrami operator 
on the projective plane}\label{multiplicity}

It was proven by the first author
in the paper~\cite{Nadirashvili1987} that the following
upper bound for the multiplicities of the eigenvalues of 
the Laplace-Beltrami operator on the projective
plane holds,
$$
m(\mathbb{RP}^2,g,\lambda_i)\leqslant 2i+3.
$$
For the first eigenvalue this means
$$
m(\mathbb{RP}^2,g,\lambda_1)\leqslant 5,
$$
which is a sharp inequality and was proved first by  
Besson~\cite{Besson1980}. 

For the second eigenvalue we have
$$
m(\mathbb{RP}^2,g,\lambda_2)\leqslant 7.
$$
The main goal of this Section is to improve the last upper bound.
\begin{Proposition}\label{multiplicity-proposition}
The following upper bound for the multiplicity of the second
eigenvalue of 
the Laplace-Beltrami operator on the projective
plane holds,
\begin{equation}\label{multiplicitybound}
m(\mathbb{RP}^2,g,\lambda_2)\leqslant 6.
\end{equation}
\end{Proposition}

For the purposes of the present paper the upper
bound~\eqref{multiplicitybound}
is sufficient. However, this bound is further improved
and generalized in the paper~\cite{Berdnikov-Nadirashvili2016}.

Let us postpone the proof and start with several lemmas.

\begin{Lemma}\label{6functions}
Let $u_1,\dots,u_6$ be linearly independent eigenfunctions
corresponding to the second eigenvalue
$\lambda_2(\mathbb{RP}^2,g).$ Then
for any point $x_0\in\mathbb{RP}^2$ there exists a 
non-trivial linear combination $v=\sum\limits_{i=1}^6\alpha_i u_i$
such that the eigenfunction $v$ has a zero of order at least $3$
at the point~$x_0.$
\end{Lemma}

\noindent{\bf Proof.} This lemma is a particular case
of Lemma 4 from paper~\cite{Nadirashvili1987}.
In fact, the proof is an easy corollary of Bers Theorem~\ref{BersTheorem}
and formula~\eqref{harmonic_pol}. $\Box$

\begin{Lemma}\label{Lemma_order}
Let $u$ be an eigenfunction
corresponding to the second eigenvalue
$\lambda_2(\mathbb{RP}^2,g)$ such that
at a point $x_1$ this eigenfunction has a zero of order
at least $3.$ Then $x_1$ is the only zero of $u$
of order greater than $1$ and the order of zero
at $x_1$ is exactly $3.$ 
\end{Lemma}

\noindent{\bf Proof.}
Since $i=2,$ $\chi(\mathbb{RP}^2)=1,$
inequality~\eqref{nodal-inequality-1}
implies in this case the inequality
$$
2\geqslant \sum\limits_{j=1}^n(m_j-1).
$$
Since $m_1\geqslant3$ and $m_i\geqslant2$ for $i>1,$  
we have $m_1-1\geqslant2,$ $m_i-1\geqslant 1$ for $i>1.$
It follows that $m_1=3$ and $n=1.$
$\Box$

Let us fix a point $x_0\in \mathbb{RP}^2$
and consider the space $V$ of eigenfunctions of $\Delta$
corresponding to the second eigenvalue
$\lambda_2(\mathbb{RP}^2,g)$ with zero of order at least $3$ at $x_0.$
Let us suppose that $\dim V\geqslant2.$ Then there exist two linearly independent
eigenfunctions $u_1,$ $u_2\in V.$ Consider
the family of eigenfunctions
\begin{equation}\label{eigenfamily}
v^\tau=\cos\tau u_1+\sin\tau u_2.
\end{equation}

\begin{Lemma}\label{unique_star}
The star $S_{x_0}(N^1(v^\tau))$ defines the eigenfunction
$v^\tau$ from formula~\eqref{eigenfamily} completely up to a sign,
i.e. if $S_{x_0}(N^1(v^{\tau_1}))=S_{x_0}(N^1(v^{\tau_2}))$ 
then $v^{\tau_1}=\pm v^{\tau_2}.$
\end{Lemma}

\noindent{\bf Proof.}
Since $S_{x_0}(N^1(v^{\tau_1}))=S_{x_0}(N^1(v^{\tau_2})),$ 
the homogeneous harmonic polynomials $P^{\tau_1}_3$ and $P^{\tau_2}_3$
corresponding by Bers Theorem~\ref{BersTheorem}
to $v^{\tau_1}$
and $v^{\tau_2}$ are proportional. But then formula~\eqref{eigenfamily}
implies that either $P^{\tau_1}_3=P^{\tau_2}_3$ or $P^{\tau_1}_3=-P^{\tau_2}_3.$
In the first case we have
$$
v^{\tau_1}-v^{\tau_2}=O(|x|^4).
$$  
Then $v^{\tau_1}-v^{\tau_2}$ is an 
eigenfunction of $\Delta$ corresponding to the second eigenvalue
$\lambda_2(\mathbb{RP}^2,g)$ with zero of order at least $4$ at $x_0.$
It follows from Lemma~\ref{Lemma_order} that $v^{\tau_1}-v^{\tau_2}\equiv0$
and therefore $v^{\tau_1}=v^{\tau_2}.$
A similar argument shows that in the second case we have $v^{\tau_1}=-v^{\tau_2}.$
$\Box$

\begin{Lemma}\label{dimension1}
Let $x_0\in\mathbb{RP}^2$ and $V$ be the space of eigenfunctions of $\Delta$
corresponding to the second eigenvalue
$\lambda_2(\mathbb{RP}^2,g)$ with a zero of order at least $3$ at $x_0.$
Then $\dim V\leqslant1.$
\end{Lemma}

\noindent{\bf Proof.} Let us suppose that $\dim V\geqslant2.$ Then
there exist two linearly independent
eigenfunctions $u_1,$ $u_2\in V.$ Consider the family of eigenfunctions
$v^\tau\in V$ defined by equation~\eqref{eigenfamily} and the family of nodal
graphs $N^1(v^\tau).$

Let $p:\mathbb{S}^2\longrightarrow\mathbb{RP}^2$ be the standard projection.
Let us consider the eigenfunction $u_1\circ p$ on the sphere $\mathbb{S}^2.$
It follows from the above mentioned arguments that the nodal graph
$N^1(u_1\circ p)$ on the sphere has the following properties:
\begin{itemize}
\item there are exactly two vertices $p^{-1}(x_0)$ that we call $N$ and $S$,
they are antipodal,
\item locally exactly $6$ edges emanate from each vertex.
\end{itemize}

\begin{Claim}\label{claim1}
All nodal domains are topological disks.
\end{Claim}

Indeed, the Euler inequality~\eqref{euler-inequality} for the nodal graph of $u_1$
on $\mathbb{RP}^2$ implies that there is at least $3$ nodal domains. In the
same time, the Courant Nodal Domain Theorem~\ref{CNDT} implies that
there are at most $3$ nodal domains. As a result, there are
exactly $3$ nodal domains
for $u_1$ on $\mathbb{RP}^2.$ Now remark that it follows that the 
Euler inequality~\eqref{euler-inequality} turns into an equality. 
It is possible if and only
if all nodal domains of $u_1$ are topological disks. Let us consider now
the nodal graph of $u_1\circ p$ on the sphere $\mathbb{S}^2.$ Since there are 
no non-trivial coverings of a disk, and the nodal domains 
of $u_1\circ p$ are preimages
of the nodal domains of $u_1,$ there are exactly $6$ nodal domains of 
$u_1\circ p$ on $\mathbb{S}^2$ and all are topological disks.
 
 \begin{Claim}\label{claim2}
The nodal graph of $u_1\circ p$ is invariant under rotation
by $\pm\frac{\pi}{3}$ around the axis going though $N$ and $S.$
\end{Claim}

Let us emphasize that ``invariant'' here and below means
``invariant up to a homotopy preserving tangent vectors at the point
$N$ and $S$''.

The proof of the Claim~\ref{claim2} is as follows.
Since $v^0=-v^{\pi},$ the nodal graph $N^1(v^\tau)$
is deformed continuously when $\tau$ changes from $0$ to $\pi$
and the result coincides with the initial graph,
$N^1(v^0)=N^1(v^{\pi}).$

Since $N^1(v^0)=N^1(v^{\pi}),$
when $\tau$ changes from $0$ to $\pi$ the $6$-ray star $S_{x_0}(N^1(v^\tau))$
rotates by
angle $k\frac{\pi}{3}.$ But then $k=\pm1.$ Indeed, if $k\ne\pm1$
then there exists $0<\tau_0<\pi$ such that
$S_{x_0}(N^1(v^{\tau_0}))$ is obtained from $S_{x_0}(N^1(v^0))$
by the rotation by angle $(\sgn k)\frac{\pi}{3}.$
Then $S_{x_0}(N^1(v^{\tau_0}))=S_{x_0}(N^1(v^0))$
and Lemma~\ref{unique_star} implies that $v^{\tau_0}=\pm v^0,$
but this contradicts the inequality $0<\tau_0<\pi.$

Let us change the direction of counting the angle in such a way that
the angle of rotation is $\frac{\pi}{3}.$ Then we have the following
result: {\em when $\tau$ changes from $0$ to $\pi,$ the star $S_{x_0}(N^1(v^\tau))$
rotates exactly by $\frac{\pi}{3}.$}

\begin{Claim}\label{claim3}
There are no loops in the nodal graph, i.e. all edges join $N$ and $S.$
\end{Claim}

Indeed, let us consider an edge $\gamma$ emanating from $N$ such that another 
endpoint of $\gamma$ is also $N.$ 
Let us numerate the vectors from the star $S_N(N^1(u_1\circ p))$
in consecutive order as $v_0,\dots,v_5$ in such a way that 
the edge $\gamma$ emanates from $N$ with the tangent 
vector $v_0.$ Then there is two cases.

In case I the tangent vector at the endpoint $N$ of $\gamma$ is 
$-v_1.$ In this case  the nodal graph is clearly not invariant under the 
rotation by $\frac{\pi}{3}.$  

Remark that the tangent
vector $-v_5$ at the endpoint $N$  could be considered as
$-v_1$ with another numeration order of the vectors from 
the $S_N(N^1(u_1\circ p)).$

In case II the tangent vector at the endpoint $N$ of $\gamma$ is $-v_k,$ 
where $k=2$ or $k=3.$ Since the nodal graph is invariant under the 
rotation by $\frac{\pi}{3},$
the edge emanating from $N$ with 
tangent vector $v_1$ has $-v_{k+1}$ as its tangent
vector at its endpoint $N.$
This implies that there are two loops
on a sphere intersecting transversally at exactly one point $N$
but this is impossible.

Remark that the tangent vector $-v_4$ at the endpoint $N$  could be considered as
$-v_2$ with another numeration order of the vectors from 
$S_N(N^1(u_1\circ p)).$

In both cases we obtain a contradiction with the assumption
that an edge can start and end at the same vertex.
Hence, all edges join $N$ to~$S.$

Let us now finish the proof of Lemma~\ref{dimension1}.
Consider the nodal graph of $u_1\circ p$ on $\mathbb{S}^2.$
A small neighbourhood of $N$ is divided by the graph in $6$ sectors,
where the signs of $u_1\circ p$ alternate.
By Claim~\ref{claim3}, these sectors lie in different nodal domains.
Since there are exactly $6$ nodal domains, there are three of them where
$u_1\circ p$ is positive and three of them where $u_1\circ p$ is negative.

Let us consider the action of the antipodal map $\sigma$ on nodal domains.
It is well-defined. Indeed, suppose $x$ and $y$ belong to the same nodal domain.
Then one can join $x$ and $y$ by a path inside their nodal domain.
Applying $\sigma$ we obtain a path joining $\sigma(x)$ and $\sigma(y),$
on which $u_1\circ p$ does not change sign. Thus, $\sigma(x)$ and $\sigma(y)$
belong to the same nodal domain.

Since $u_1\circ p$ is obtained from the eigenfunction $u_1$ on $\mathbb{RP}^2,$
the antipodal map preserves the sign of $u_1\circ p.$ Since there are three nodal domains
of the same sign, there is at least one nodal domain that is mapped
by $\sigma$ to itself. At the same time, by Claim~\ref{claim1} each nodal domain
is a topological disk. Since $\sigma$ is a free involution, it can not map a disk
into itself by Brouwer's theorem. This completes the proof of Lemma~\ref{dimension1}.
$\Box$

\noindent{\bf Proof} of Proposition~\ref{multiplicity-proposition}.
Let us suppose that $m(\mathbb{RP}^2,g,\lambda_2)>6.$
Then there exist 7 linearly independent eigenfunctions
$\varphi_1,\dots,\varphi_7$ corresponding to the second eigenvalue
$\lambda_2(\mathbb{RP}^2,g).$

Let us fix a point $x_0\in\mathbb{RP}^2.$ Let us apply Lemma~\ref{6functions}
to $\varphi_1,\dots,\varphi_6$ and obtain an eigenfunction 
$u_1=\sum\limits_{i=1}^6\alpha_i\varphi_i$
with zero of order at least $3$ at the point $x_0.$ Then by Lemma~\ref{Lemma_order}
the point $x_0$ is a zero of order exactly $3.$

We can suppose without loss of generality that $\alpha_1\ne0.$
Let us then apply Lemma~\ref{6functions}
to the eigenfunctions $\varphi_2,\dots,\varphi_7$ and obtain an eigenfunction 
$u_2=\sum\limits_{i=2}^7\beta_i\varphi_i$
with zero of order at least $3$ at the point $x_0.$ Then by Lemma~\ref{Lemma_order}
the point $x_0$ is a zero of order exactly $3.$ 

Let us remark that $u_1$
and $u_2$ are linearly independent since $\alpha_1\ne0.$ 
This contradicts Lemma~\ref{dimension1}.
$\Box$

\section{Harmonic maps with branch points and metrics with
conical singularities}\label{harmonic-conical}

Let us recall the definition of a harmonic map,
see e.g. the review~\cite{Eells-Lemaire1978}.

\begin{Definition}
Let $(M,g)$ and $(N,h)$ be Riemannian manifolds. A smooth
map $f:M\longrightarrow N$ is called
harmonic if $f$ is an extremal for the energy functional
\begin{equation}\label{energy-functional}
E[f]=\int_M|df(x)|^2\,dVol_g.
\end{equation}
\end{Definition}

The following theorem (see, e.g. the paper~\cite{Eells-Sampson1964})
explains the relation between minimal and harmonic maps
in the class of isometric immersions.

\begin{Theorem}\label{harmonic-minimal}
Let $M,$ $N$ be Riemannian manifolds. If $f:M\looparrowright N$
is an isometric immersion, then $f$ is harmonic if and only
if $f$is minimal.
\end{Theorem}

Theorem~\ref{Nadirashvili-ElSoufi-Ilias} and 
Theorem~\ref{harmonic-minimal} imply the following Proposition.

\begin{Proposition}\label{metrics-harmonic-maps}
The extremal metrics on a compact surface $M$ are exactly the
metrics induced on $M$ by harmonic
immersions $M\looparrowright\mathbb{S}^n.$
\end{Proposition}

It turns out however that it is useful to consider
a wider class of harmonic immersions with branch points.

\begin{Definition}[see e.g.~\cite{Gulliver-Osserman-Royden1973}]
Let $M$ be a manifold of dimension $2.$ A~smooth map $f:M\longrightarrow N$
has a branch point of order $k$ at point $p$
if there exist local coordinates $u_1,$ $u_2$ centered at $p$ and defined in a
neighborhood of $p$ and $x_1,\dots,x_n$ centered at $f(p)$ and defined
in a neighborhood of $f(p)$ such that in these coordinates $f$
is written as
\begin{gather*}
x_1+ix_2=w^{k+1}+\sigma(w),\\
x_k=\chi_k(w),\quad k=3,\dots,n,\\
\sigma(w),\chi_k(w)=o(|w|^{k+1}),\\
\frac{\partial\sigma}{\partial u_j}(w),
\frac{\partial\chi_k}{\partial u_j}(w)=o(|w|^k),\quad j=1,2,\quad k=3,\dots,n,
\end{gather*}
where $w=u^1+iu^2.$
\end{Definition}

If $M$ is compact then a map $f:M\longrightarrow N$ could have only
finite number of branch points.

However we have now a problem. If $(N,g)$ is a Riemannian manifold
and $f:M\looparrowright N$ is an immersion with branch points, then
the induced metric $f^*g$ is not a smooth metric.

\begin{Definition}[see e.g.~\cite{Kokarev2014}]
A point $p$ on a surface $M$ is called a conical singularity 
of order $\alpha>-1$ or angle $2\pi(\alpha+1)$
of a metric $g$ if in an appropriate local complex coordinate
$z$ centered at $p$ the metric has the form
$$
g(z)=|z|^{2\alpha}\rho(z)|dz|^2
$$
in a neighborhood of $p,$ where $\rho(0)>0.$
\end{Definition}

Then we obtain immediately the following Proposition.
\begin{Proposition}
If $M$ is a compact surface, $(N,h)$ is a Riemannian manifold
and $f:M\looparrowright N$ is an immersion with branch point, 
then $g=f^*h$ is a smooth Riemannian metric except a finite
number of branch points of the map $f.$ At these points the
metric $g$ has conical singularities. The order of the conical
singularity at a point $p$ is equal to the order of 
$p$ as a branch point.
\end{Proposition}

Thus, we switch to a setting larger than the initial one.
We consider not only Riemannian metrics but also Riemannian
metrics with a finite number of conical singularities and
not only harmonic immersions but also harmonic immersions
with branch points. Then we should check that all
key results from the previous sections hold. 

It is well-known that the eigenvalues of the Laplace-Beltrami
operator could be defined using a variational approach,
\begin{equation}\label{Rayleigh-variational-formula}
\lambda_k=\min_{\stackrel{V\subset H^1(M)}{\dim V=k}}%
\max_{\stackrel{u\in V}{u\perp 1}}R[v],
\end{equation}
where
$$
R[v]=\frac{\int_M|\nabla u|^2\,dVol}{\int_M|u|^2\,dVol}
$$
is the Rayleigh quotient. This formula holds also
in the case of metrics with conical singularities,
see e.g.~\cite{Kokarev2014}.

\begin{Proposition}[\protect{\cite[Corollary 4.7]{Kokarev2014a}}]
Theorem~\ref{Nadirashvili-ElSoufi-Ilias}
holds if we consider metrics with conical singularities
and harmonic maps with branch points.
\end{Proposition}

The next problem is to prove that $V,$ $E$ and
$F$ are finite and inequality~\eqref{euler-inequality}
holds. The problem is that in the case of surfaces with isolated 
conical singularities the points of $N^2(u)$ can a priori
accumulate towards singularities. It turns out that it is not possible
since this possibility can be ruled out using resolution procedure
used in the 
papers~\cite[Lemma 3.1.1]{Karpukhin-Kokarev-Polterovich2014}
and~\cite{Kokarev2014}
in order to prove the finiteness of a nodal graph in other
contexts.

Let us define the resolution procedure 
following the paper~\cite{Karpukhin-Kokarev-Polterovich2014}. 
Let $x\in N^2(u)$ be a vertex of nodal graph. If $n=\ord_x(u)$ then 
the degree of this vertex is $2n.$
According to Bers's Theorem~\ref{BersTheorem} there exists
a neighborhood $U$  of $x$ diffeomorphic to a disk such
that $U$ does not contain other vertices and such that nodal arcs incident to 
$x$ intersect $U$ at $2n$ points precisely.
Let us denote these intersection points by $y_i,$ where $i = 0,\dots,2n-1,$ 
and assume that they are ordered consequently in the clockwise fashion. 
A new graph is obtained from the nodal graph by changing it inside $U$ 
and removing possibly appeared edges without vertices. More precisely, we 
remove the nodal set inside $U$ and round-off the edges
on the boundary $U$ by non-intersecting arcs in $U$ joining the points 
$y_{2j}$ and $y_{2j+1}.$
If there was an edge that starts and ends at $x,$ then such a procedure may 
make it into a loop. If this occurs, then we remove this loop to obtain a 
genuine graph in the sufrace.
The new graph has one vertex less and at most as many faces as the original graph.

We give now a short proof by Karpukhin (with his permission).

\begin{Proposition}[Karpukhin,~\cite{Karpukhin2013q}]
A nodal graph of an eigenfunction $u$ on a surface $M$ with isolated 
conical singularities is finite.
\end{Proposition}

\noindent{\bf Proof.} Suppose that there are infinitely many points in 
$N^2(u)$, it is easy to see that in this case the set $N^2(u)$ is countable. 
Then the only possible accumulation points of $N^2(u)$ are conical 
singularities. For each conical singularity $p_j$ let us choose a base 
of neighbourhoods $V_i^{(j)}$ such that $\bar V_{i+1}^{(j)}\subset V_i^{(j)}$ 
and $N^2(u)\cap\bigcup\limits_{i=1}^\infty\partial V_i^{(j)}=\emptyset$. 
Hence for the sets $V_i=\bigcup\limits_j V_i^{(j)}$ we have
$\bar{V}_{i+1}\subset V_{i}$, $N^2(u)\cap\bigcup\limits_{i=1}^\infty \partial V_i%
=\emptyset$ and $M\backslash V_i$ contains only finite quantity 
of elements of $N^2(u)$. For any $i$ for the points of $N^2(u)$ in 
$V_i\backslash\bar{V}_{i+1}$ one can choose a collection of 
disjoint neighbourhoods $U_{ki}$ such that 
$\bar{U}_{ki}\subset V_i\backslash\bar{V}_{i+1}$. 
Thus we constructed a collection of disjoint neighbourhoods of all points 
in $N^2(u)$.     

Next we apply the resolution procedure at all but finite number of vertices. 
Choosing this finite number big enough and 
applying Euler's inequality we arrive at contradiction with Courant's 
nodal domain theorem. $\Box$

Thus, in the setting of metrics with conical singularities
inequality~\eqref{euler-inequality} and all results
obtained with its help hold, including the key upper
bound $m(\mathbb{RP}^2,g,\lambda_2)\leqslant 6$
from Proposition~\ref{multiplicity-proposition}
from Section~\ref{multiplicity}.

Let us also remark that for any manifold equipped with a metric
with isolated conical singularities 
it is possible to construct
a sequence of smooth Riemannian manifolds such that 
their area as well as their eigenvalues
converge to the area and eigenvalues of the initial manifold,
see e.g.~\cite{Sher2015}.

\section{Calabi-Barbosa theorem and its 
implications}\label{Calabi-Barbosa-implications}

Now we should study harmonic immersions with branched
points $\mathbb{RP}^2\looparrowright\mathbb{S}^n.$
Since we  have the upper bound $m(\mathbb{RP}^2,g,\lambda_2)\leqslant 6$
from Proposition~\ref{multiplicity-proposition},
all immersions corresponding to $\lambda_2$ are among
immersions $\mathbb{RP}^2\looparrowright\mathbb{S}^5.$

Let $p:\mathbb{S}^2\longrightarrow\mathbb{RP}^2$ be the standard projection.
We can lift a harmonic immersion with branch points 
$f:\mathbb{RP}^2\looparrowright\mathbb{S}^5$
to a harmonic immersion with branch points 
$F=f\circ p:\mathbb{S}^2\longrightarrow\mathbb{S}^5.$

The following theorem was proved by Calabi in 1967 and later refined
by Barbosa in 1975. 
Let $g_{\mathbb{S}^n}$ denote the standard metric on $\mathbb{S}^n.$
The radius of $\mathbb{S}^n$ is $1.$

\begin{Theorem}[Calabi~\cite{Calabi1967}, 
Barbosa~\cite{Barbosa1975}]\label{Calabi-Barbosa}
Let $F:\mathbb{S}^2\longrightarrow\mathbb{S}^n$
be a harmonic immersion with branch points such
that the image is not contained in a hyperplane.
Then
\begin{itemize}
\item[(i)] the area of\/ $\mathbb{S}^2$ with
respect to the induced metric
$\Area(\mathbb{S}^2,F^*g_{\mathbb{S}^n})$
is an integer multiple of $4\pi;$
\item[(ii)] $n$ is even, $n=2m,$ and
$$
\Area(\mathbb{S}^2,F^*g_{\mathbb{S}^n})\geqslant2\pi m(m+1).
$$
\end{itemize}
\end{Theorem}

\begin{Definition}
If $\Area(\mathbb{S}^2,F^*g_{\mathbb{S}^n})=4\pi d,$ then we say that $F$
is of harmonic degree $d.$
\end{Definition}

We obtain immediately a lower bound for the harmonic degree.

\begin{Proposition}
Let $F:\mathbb{S}^2\longrightarrow\mathbb{S}^{2m}$
be a harmonic immersion with branch points such
that the image is not contained in a hyperplane.
Then $d\geqslant\frac{m(m+1)}{2}.$
\end{Proposition}

Calabi-Barbosa Theorem~\ref{Calabi-Barbosa} implies
the following Proposition.

\begin{Proposition}\label{sufficient-maps}
It is sufficient for our goals to consider harmonic immersions
with branch points $F:\mathbb{S}^2\longrightarrow\mathbb{S}^4$
(such that the image is not contained in a hyperplane)
of harmonic degree $d\geqslant3$ 
and  $F:\mathbb{S}^2\longrightarrow\mathbb{S}^2.$ 
\end{Proposition}

It follows that we should consider only harmonic immersions
with branch points $\mathbb{RP}^2\longrightarrow\mathbb{S}^2$
and $\mathbb{RP}^2\longrightarrow\mathbb{S}^4.$ However, the following
Proposition permits to exclude maps $\mathbb{RP}^2\longrightarrow\mathbb{S}^2.$

\begin{Proposition}[see e.g.~\cite{Eells-Lemaire1978}]\label{no-rp2-s2}
Every harmonic map
$\mathbb{RP}^2\longrightarrow\mathbb{S}^2$
is constant.
\end{Proposition}

\section{Harmonic maps from $\mathbb{S}^2$ to $\mathbb{S}^4$
and their singularities}\label{space_of_harmonic_maps}

Let us recall the well-known Penrose twistor map
$$
T:\mathbb{CP}^3\longrightarrow\mathbb{HP}^1\cong\mathbb{S}^4,\quad
T([z_0:z_1:z_2:z_3])=[z_0+z_1j:z_2+z_3j].
$$
Let $z$ be a conformal parameter on $\mathbb{S}^2.$
\begin{Definition}
Let us call a curve 
$$
f:\mathbb{S}^2\longrightarrow\mathbb{CP}^3,\quad f(z)=[f_0(z):f_1(z):f_2(z):f_3(z)],
$$
horizontal if
$$
f'_1f_2-f_1f'_2+f'_3f_4-f_3f'_4=0.
$$
\end{Definition}

In 1982, Bryant described in the paper~\cite{Bryant1982} 
a very important relation between harmonic immersions
with branch points $\mathbb{S}^2\longrightarrow\mathbb{S}^4$
and (anti)holomorphic horizontal curves in $\mathbb{CP}^3.$

Let $A:\mathbb{S}^4\longrightarrow\mathbb{S}^4$ be the antipodal map.

\begin{Theorem}[Bryant~\cite{Bryant1982}]\label{Bryant-theorem}
For each harmonic immersion
with branch points $F:\mathbb{S}^2\longrightarrow\mathbb{S}^4$
there exists either a holomorphic or an antiholomorphic
horizontal curve $f:\mathbb{S}^2\longrightarrow\mathbb{CP}^3,$
such that $T\circ f=F,$
$$
\xymatrix{
&\mathbb{CP}^3\ar[d]^T\\
\mathbb{S}^2\ar[r]^F\ar[ur]^f & \mathbb{S}^4
}
$$
For each (anti)holomorphic
horizontal curve $f:\mathbb{S}^2\longrightarrow\mathbb{CP}^3$
the map $F=T\circ f:\mathbb{S}^2\longrightarrow\mathbb{S}^4$
is a harmonic immersion with branch points.

If a harmonic immersion $F:\mathbb{S}^2\longrightarrow\mathbb{S}^4$
has a holomorphic (antiholomorphic)
horizontal curve $f:\mathbb{S}^2\longrightarrow\mathbb{CP}^3,$
then $A\circ F:\mathbb{S}^2\longrightarrow\mathbb{S}^4$
has an antiholomorphic (holomorphic) horizontal curve.
\end{Theorem}

\begin{Definition}
An (anti)holomorphic horizontal curve $f$ appearing in Bryant's
Theorem~\ref{Bryant-theorem} is called
the lift of an harmonic immersion $F.$
\end{Definition}

Let us remark that $F$ and $A\circ F$ induce the same metric
on $\mathbb{S}^2.$ It follows that it is sufficient to consider
harmonic immersions with holomorphic lifts.

\begin{Theorem}[Bryant~\cite{Bryant1982}]\label{Bryant-theorem-degree}
Let $F:\mathbb{S}^2\longrightarrow\mathbb{S}^4$ be a harmonic
immersion with branched points of harmonic degree $d$
with holomorphic lift $f:\mathbb{S}^2\longrightarrow\mathbb{CP}^3.$
Then $f:\mathbb{S}^2\longrightarrow\mathbb{CP}^3$ is an algebraic
curve of degree $d.$
\end{Theorem}

Now we need some results from the theory of higher singularities
of these holomorphic horizontal lifts, see
e.g. the paper~\cite{Bolton-Woodward2000} by Bolton
and Woodward. Let $[\mathbf{f}(z)]=[f_0(z),f_1(z),f_2(z),f_3(z)]$
be a representative of $f:\mathbb{S}^2\longrightarrow\mathbb{CP}^3$
in the homogeneous coordinates in $\mathbb{CP}^3.$
Let $\mathbf{f}^{(i)}(z)$ denotes the $i$th derivative of $\mathbf{f}(z).$
Let 
$$
Z(f)=\{z\,|\,%
\mathbf{f}(z)\wedge\mathbf{f}'(z)\wedge\ldots\wedge\mathbf{f}^{(3)}(z)=0\}.
$$
Remark that $Z(f)$ consists of isolated points if $f$ is linearly full,
i.e. if the image of $f$ is not inside a hyperplane.

Let us apply the Gram-Schmidt orthogonalization process
to $\mathbf{f}(z),$ $\mathbf{f}'(z),$ $\mathbf{f}''(z),$ $\mathbf{f}'''(z)$
at $z\not\in Z(f)$ and obtain
$\tilde{\mathbf{f}}_0(z)=\mathbf{f}(z),$ $\tilde{\mathbf{f}}_1(z),$
$\tilde{\mathbf{f}}_2(z),$ $\tilde{\mathbf{f}}_3(z).$
Then it turns out that the trivial bundle $\mathbb{S}^2\times\mathbb{C}^4$
has an orthogonal decomposition as a sum of holomorphic linear bundles
$$
\mathbb{S}^2\times\mathbb{C}^4=L_0\oplus\ldots\oplus L_3,
$$
such that $L_i$ is spanned by $\tilde{\mathbf{f}}_i$
for $z\not\in Z(f).$ These $L_i$ describe the Fr\'enet frame for $f.$

The bundle map 
$\partial_i:T^{1,0}\mathbb{S}^2\otimes L_i\longrightarrow L_i^\perp$
given by
$$
\partial_i\left(\frac{\partial}{\partial z}\otimes s_i\right)=%
\left(\frac{\partial s_i}{\partial z}\right)^\perp,
$$
where $s_i$ is a local holomorphic section of $L_i,$ and 
$\left(\dfrac{\partial s_i}{\partial z}\right)^\perp$
denotes the component of $\dfrac{\partial s_i}{\partial z}$
orthogonal to $L_i,$
satisfies
$$
\partial_i\left(\dfrac{\partial}{\partial z}\otimes\tilde{\mathbf{f}}_i(z)\right)=%
\tilde{\mathbf{f}}_{i+1}(z).
$$
It follows that $\partial_i$ is a holomorphic map and has the image in $L_{i+1}.$

\begin{Definition}
A (linearly full) holomorphic curve
$f:\mathbb{S}^2\longrightarrow\mathbb{CP}^3$ has a higher singularity
of type $(r_0(p),r_1(p),r_2(p))$
at a point $p\in Z(f)$ if for $i=0,1,2$ the holomorphic bundle
maps
$\partial_i$ has a zero of order $r_i(p)$
at $p$ and $r_0(p)+r_1(p)+r_2(p)>0.$
\end{Definition}

It turns out that for a horizontal
curve one has $r_2(p)=r_0(p),$ i.e. its higher singularity type
at a point $p$ is described by two integers $r_0(p),$ and 
$r_1(p).$ Let us define quantities
$$
r_0=\sum_{p}r_0(p),\quad r_1=\sum_{p}r_1(p).
$$
The next Proposition relates them to the degree $d.$

\begin{Proposition}[Bolton, 
Woodward~\cite{Bolton-Woodward2001}]\label{singularities-proposition}
For a linearly full holomorphic horisontal curve in $\mathbb{CP}^3$
the following equation holds,
$$
2r_0+r_1=2d-6.
$$
\end{Proposition}
 
We need here to recall the definition
of an umbilic point.

\begin{Definition}\label{umbilic-definition}
Let $(M,g)$ and $(N,h)$ be Riemannian manifolds
and $\nabla^M$ and $\nabla^N$ be the corresponding
Levi-Civita connections. 

Let $F:M\longrightarrow N$ be an immersion. Then a) the
second fundamental form $\mathbf{II}^F$ of $F$ is defined by the
formula
$$
\nabla^N_{dF(X)}dF(Y)=dF(\nabla^M_XY)+\mathbf{II}^F(X,Y);
$$
b) the vector field 
$$
\zeta=\frac{1}{\dim M}\tr\mathbf{II}^F
$$
is called a mean curvature normal vector;

c) a point $p\in M$ is called an umbilic point if there exists a vector $v\in T_{F(p)}N$
such that at the point $p$ one has
\begin{equation}\label{umbilic-1}
\mathbf{II}^F_p(X,Y)=g_p(X,Y)\cdot v.
\end{equation}
\end{Definition}

It follows immediately from Definition~\ref{umbilic-definition}
that if $p$ is an umbilic
then $\mathbf{II}^F_p(X,Y)=g_p(X,Y)\cdot\zeta(p).$

As an example it is useful to consider a classical case of
an immersion $F$ of a two-dimensional surface $M$ 
to $N=\mathbb{R}^3$ equipped with the euclidean metric $h.$
Let us consider the induced metric $g=F^*h$ on $M.$
Then it is easy to check that 
$\mathbf{II}^F(X,Y)=II(X,Y)\cdot\vec{n},$
where $II(X,Y)$ is the classical second fundamental form
of the surface $M$ and $\vec{n}$ is a unit normal vector
field on $M.$ Let us recall that in the basis consisting
of principal directions the metric $g$ has the identity matrix
and the classical second fundamental form $II$ has
the diagonal matrix with the principal curvatures
$\lambda_1$ and $\lambda_2$ on the diagonal.
Then formula~\eqref{umbilic-1} is equivalent to the
equality $\lambda_1=\lambda_2$ which is the classical
definition of an umbilic point for a two-dimensional surface
in the Euclidean space $\mathbb{R}^3.$

Let $z$ be a conformal parameter on $\mathbb{S}^2.$
It is easy to check that the following Proposition holds.

\begin{Proposition}
A point $p\in\mathbb{S}^2$ is an umbilic point of
a harmonic immersion $F:\mathbb{S}^2\longrightarrow\mathbb{S}^4$
if and only if 
\begin{equation}\label{II-zero}
\mathbf{II}^F_p(\partial/\partial z,\partial/\partial z)=0.
\end{equation}
\end{Proposition}

\noindent{\bf Proof.} If a point $p$ is umbilic then
$$
\mathbf{II}^F_p(\partial/\partial z,\partial/\partial z)=%
g_p(\partial/\partial z,\partial/\partial z)\cdot\zeta(p)=0,
$$
since $z$ is a conformal coordinate and $g_p=2\Phi|dz|^2$
for some $\Phi.$

Let equality~\eqref{II-zero} holds. Since $F$ is real, this implies that
$$
\mathbf{II}^F_p(\partial/\partial\bar{z},\partial/\partial\bar{z})=0.
$$
If follows that formula~\eqref{umbilic-1} holds for
$$
v=%
\frac{\mathbf{II}^F_p(\partial/\partial z,\partial/\partial\bar{z})}%
{g_p(\partial/\partial z,\partial/\partial\bar{z})},
$$
and $p$ is umbilic.
$\Box$

The higher singularities of a holomorphic horizontal lift $f$ of a harmonic
immersion with branched points $F:\mathbb{S}^2\longrightarrow\mathbb{S}^4$
are related to the branch points and the umbilics of $F.$

\begin{Proposition}[Bolton, 
Woodward~\cite{Bolton-Woodward2000,
Bolton-Woodward2001}]\label{singularities-relation}
A point $p$ is a branch point of $F$ if and only if $r_0(p)>0.$
Moreover, $r_0(p)$ is equal to the order of zero of 
$dF(\partial/\partial z)$ at $p.$

If $r_0(p)=0$ then $p$ is an umbilic if and only if $r_1(p)>0.$
Moreover, $r_1(p)$ is equal to the order of zero of\/
$\mathbf{II}^F(\partial/\partial z,\partial/\partial z)$ at $p.$

The higher singularities of $f$ occur exactly at the branch points
and umbilics of $F.$
\end{Proposition}

Combining Propositions~\ref{singularities-proposition}
and~\ref{singularities-relation}, we obtain the following Proposition.

\begin{Proposition}\label{s2-s4}
Let $F:\mathbb{S}^2\longrightarrow\mathbb{S}^4$ be a harmonic
immersion with branch points of harmonic degree $d.$ Then
\begin{itemize}
\item[(i)] if $d=3$ then $F$ does not have either branch points or
umbilics,
\item[(ii)] if $d>3$ then $F$ has at least one branch point or an umbilic.
\end{itemize}
\end{Proposition}

\section{Existence of maximal metrics}\label{existence}

What can we say about the existence of the maximal metric
for a given eigenvalue on a given surface? The situation
in the case of the first eigenvalue is the following.

\begin{Theorem}[Matthiesen, Siffert~\cite{Matthiesen-Siffert2017q}]
For any closed surface $M,$ there is a metric $g$ on $M,$
smooth away from finitely many conical singularities,
achieving $\Lambda_1(M),$ i.e.
$$
\Lambda_1(M)=\bar{\lambda}_i(M,g)=\lambda_1(M,g)\Area(M,g).
$$
\end{Theorem}

However, as we observed in the Introduction, the situation is
more complicated for higher eigenvalues. In particular,
on the sphere there is no maximal metrics for $\bar{\lambda}_k$
if $k>1,$ see the papers~\cite{Nadirashvili2002,Petrides2014} for $k=2,$
\cite{Nadirashvili-Sire2015q2} for $k=3$ and
\cite{Karpukhin-Nadirashvili-Penskoi-Polterovich2017q}
for arbitrary $k>1.$

It turns out that extremal metrics for higher eigenvalues
on the sphere exhibit the so-called ``bubbling phenomenon''. 
This phenomenon was studied in details by the first author and Sire in the 
papers~\cite{Nadirashvili-Sire2015a,Nadirashvili-Sire2015b}
and also by Petrides~\cite{Petrides2017}
in the context of maximization of eigenvalues in a given
conformal class. More precisely,  
they investigated the question of existence of 
Riemannian metrics with conical singularities for 
which the quantity 
$$
\Lambda_k(M,[g]) = \sup_{h\in [g]}\bar\lambda_k(M,h)
$$
is attained, where $[g]$ denotes the class of metrics conformally 
equivalent to $g.$

The equality
$$
\Lambda_k(\mathbb{S}^2)=8\pi k.
$$
proven in the recent 
paper~\cite{Karpukhin-Nadirashvili-Penskoi-Polterovich2017q}
combined with~\cite[Theorem 2]{Petrides2017}
implies the following result.

\begin{Proposition}[\cite{Karpukhin-Nadirashvili-Penskoi-Polterovich2017q}]
Let $(M,g)$ be a closed Riemannian surface and $k\geqslant 2$. If 
$$
\Lambda_k(M,[g])>\Lambda_{k-1}(M,[g]) + 8\pi,
$$
then there exists a maximal metric $\tilde g\in [g]$, smooth 
except possibly at a finite set of conical singularities, such 
that $\Lambda_k(M,[g]) = \bar\lambda_k(M,\tilde g)$.
\end{Proposition}

Since there is only one conformal structure on $\mathbb{RP}^2,$
see e.g. the book \cite{Schifer-Spencer1954},
$\Lambda_2(\mathbb{RP}^2,[g])=\Lambda_2(\mathbb{RP}^2)$ and
we have the following Proposition.

\begin{Proposition}\label{rp2-max}
If 
$$
\Lambda_2(\mathbb{RP}^2)>\Lambda_1(\mathbb{RP}^2) + 8\pi=20\pi,
$$
then there exists a maximal metric $\tilde g$, smooth 
except possibly at a finite set of conical singularities, such 
that $\Lambda_2(\mathbb{RP}^2)=\bar\lambda_2(\mathbb{RP}^2,\tilde g)$.
\end{Proposition}

\section{Proof of Theorem~\ref{maintheorem}}\label{mainproof}

Let us consider a sequence $\{g_n\}$ of metrics of area one 
on the projective plane such that the limiting metric is a singular 
metric realized as a union of the projective plane and the sphere 
touching at a point, with standard metrics and the ratio of the areas~$3:2.$
In this case the limit spectrum is the union of spectra
of the projective plane with standard metric $g'$
of area $\frac{3}{5}$ and
of the sphere with standard metric $g''$ of area $\frac{2}{5},$
see e.g.~\cite[Section 2]{Colbois-ElSoufi2003} and~\cite{Anne1987} for 
more details about the limit spectrum.
Then
$$
\lim_{n\to\infty}\lambda_2(\mathbb{RP}^2,g_n)=\lambda_1(\mathbb{RP}^2,g')%
=\lambda_1(\mathbb{S}^2,g'')=20\pi.
$$
Hence,
$$
\lim_{n\to\infty}\bar{\lambda}_2(\mathbb{RP}^2,g_n)=20\pi.
$$

If $\Lambda_2(\mathbb{RP}^2)=20\pi,$ then the proof is finished.
If $\Lambda_2(\mathbb{RP}^2)>20\pi,$
then by Proposition~\ref{rp2-max}
there exists a maximal metric $\tilde g$, smooth 
except possibly at a finite set of conical singularities, such 
that 
$$
\Lambda_2(\mathbb{RP}^2)=\bar\lambda_2(\mathbb{RP}^2,\tilde g)
$$
As we already know from Proposition~\ref{metrics-harmonic-maps},
the metric $\tilde g$ is induced on $\mathbb{RP}^2$
by a harmonic immersion with branched points
$\mathbb{RP}^2\longrightarrow\mathbb{S}^n.$

The upper bound $m(\mathbb{RP}^2,g,\lambda_2)\leqslant 6$
from Proposition~\ref{multiplicity-proposition} implies that
all harmonic immersions with branch points
corresponding to $\lambda_2$ are among
immersions $\mathbb{RP}^2\looparrowright\mathbb{S}^5.$

Let $p:\mathbb{S}^2\longrightarrow\mathbb{RP}^2$ be the 
standard projection. We can lift a harmonic immersion with branch points 
$f:\mathbb{RP}^2\looparrowright\mathbb{S}^5$
to a harmonic immersion with branch points 
$F=f\circ p:\mathbb{S}^2\longrightarrow\mathbb{S}^5.$

Calabi-Barbosa's Theorem~\ref{Calabi-Barbosa} implies
that it is sufficient to consider 
a harmonic immersion with branch points $F:\mathbb{S}^2\longrightarrow\mathbb{S}^4$
of harmonic degree $d\geqslant3$ such
that the image is not contained in a hyperplane
and  a harmonic immersion $F:\mathbb{S}^2\longrightarrow\mathbb{S}^2,$
see Proposition~\ref{sufficient-maps}.
However, Proposition~\ref{no-rp2-s2}
says that we can exclude harmonic maps
$\mathbb{RP}^2\longrightarrow\mathbb{S}^2$
since they are constant.
As a result, we should consider only
a harmonic immersion with branch points
$\mathbb{RP}^2\longrightarrow\mathbb{S}^4.$

Consider a harmonic immersion with branch points
$f:\mathbb{RP}^2\longrightarrow\mathbb{S}^4$
corresponding to $\lambda_2$ and its lift
$F=f\circ p:\mathbb{S}^2\longrightarrow\mathbb{S}^4.$
As we know from Proposition~\ref{s2-s4},
there are two different cases depending on its
harmonic degree $d.$

Consider the case $d=3.$
Let $g_{\mathbb{S}^n}$ denote the standard metric on~$\mathbb{S}^n.$
Since $d=3,$ one has $\Area(\mathbb{S}^2,F^*g_{\mathbb{S}^n})=12\pi$
due to Calabi-Barbosa Theorem~\ref{Calabi-Barbosa}.
Then $\Area(\mathbb{RP}^2,f^*g_{\mathbb{S}^n})=6\pi$
because $p:\mathbb{S}^2\longrightarrow\mathbb{RP}^2$
is a two-sheeted covering. Since the radius of $\mathbb{S}^n$ is $1,$
Takahashi Theorem~\ref{takahashi} implies that $\lambda_2=2.$
As a result, $\bar{\lambda}_2(\mathbb{RP}^2,f^*g_{\mathbb{S}^n})=12\pi<20\pi$
and the induced metric is not maximal.

Consider the case $d>3.$ In this case Proposition~\ref{singularities-proposition}
implies that $F=f\circ p:\mathbb{S}^2\longrightarrow\mathbb{S}^4$
and hence $f:\mathbb{RP}^2\longrightarrow\mathbb{S}^4$ have at least one
branch point or umbilic. Let us prove that an immersion by eigenfunctions
corresponding to $\lambda_2$ cannot have either branch points or umblilics.

Let us suppose that $f=(f^1,\dots,f^5)$ and $p\in\mathbb{RP}^2$ is a branch point.
It follows that $f^i$ are linearly independent
eigenfunctions with eigenvalue $\lambda_2=2$
such that $df^i(p)=0.$ One can then construct at least $4$
linearly independent eigenfunctions $\tilde{f}^i,$ $i=1,\dots,4,$
such that $\tilde{f}^i(p)=0,$ $d\tilde{f}^i(p)=0.$
This means that all $\tilde{f}^i$ have zero of order $2$ at $p.$
Using Bers Theorem~\ref{BersTheorem} one can then construct at least $2$
linearly independent eigenfunctions with eigenvalue $\lambda_2=2$
with zero of order $3$ at $p,$ but this contradicts 
Lemma~\ref{dimension1}.

Let us suppose that $f=(f^1,\dots,f^5)$ and $p\in\mathbb{RP}^2$
is an umbilic. Let $z$ be a local
conformal parameter on $\mathbb{RP}^2$ in a neighborhood of the point $p.$
Let $ds^2=2\Phi|dz|^2$ be the induced metric.
It is well-known that $f_{z\bar{z}}=-\Phi f,$ see 
e.g.~\cite{Barbosa1975,Calabi1967}, this is in fact a harmonic map
equation in this particular setting.
Since $p$ is an umbilic, 
$\mathbf{II}^f_p(\partial/\partial z,\partial/\partial z)=0.$
By definition of the second fundamental form, this means that $f_{zz}(p)$ 
is a tangent vector and hence $f_{zz}(p)$ is a linear combination of 
$f_z(p)$ and $f_{\bar{z}}(p).$ It follows that there exist $\alpha,$ $\beta\in\mathbb{C}$
such that for any $i=1,\dots,5$ the following equations hold,
\begin{gather}
f^i_{z\bar{z}}(p)=-\Phi(p)f^i(p),\label{eq1}\\ 
f^i_{zz}(p)=\alpha f^i_z(p)+\beta f^i_{\bar{z}}(p),\label{eq2}\\
f^i_{\bar{z}\bar{z}}(p)=\bar{\beta} f^i_z(p)+%
\bar{\alpha} f^i_{\bar{z}}(p).\label{eq3}
\end{gather}
Remark that these equations are linear. This implies that they hold
for any linear combination of $f^i.$

Now one can construct two linear combinations 
$$
\varphi=\sum_{i=1}^5 A_if^i,\quad \psi=\sum_{i=1}^5 B_if^i
$$
with real coefficients $A_i,$ $B_i$ such that $\varphi$ and
$\psi$ have zero of order $2$ at~$p.$ It follows that
$$
\varphi(p)=\varphi_z(p)=\varphi_{\bar{z}}(p)=0,\quad
\psi(p)=\psi_z(p)=\psi_{\bar{z}}(p)=0.
$$
As it was remarked before, the equations
$\eqref{eq1},$ $\eqref{eq2},$ $\eqref{eq3}$ hold for 
$\varphi$ and $\psi.$ It follows that they
are eigenfunctions with eigenvalue $\lambda_2=2$ with
zero of order $3$ at the point $p.$
This contradicts Lemma~\ref{dimension1}.

Thus, it is proven that for any extremal
metric $g$ smooth except possibly finite number of  conical singularities
one has
$\bar{\lambda}_2(\mathbb{RP}^2,g)=12\pi.$
This contradicts our assumption 
$\bar{\lambda}_2(\mathbb{RP}^2,g)>20\pi$
and finishes the proof. $\Box$

\end{document}